\documentclass[12pt]{article}
\usepackage{amsfonts}

\newtheorem{preqn}{Question}[section]
\newenvironment{qn}{\preqn\rm}{\endpreqn}

\newtheorem{theorem}{Theorem}[section]
\newcommand{\Aut}{\mathop{\mathrm{Aut}}}
\newcommand{\Cay}{\mathop{\mathrm{Cay}}}

\begin{document}
\title{Around homogeneity}
\author{Peter J. Cameron\\University of St Andrews}
\date{To the memory of Robert Woodrow}
\maketitle

\begin{abstract}
Forty-five years ago, a young researcher in finite permutation group theory
encountered a paper by Robert Woodrow. The homogeneous triangle-free graph
Woodrow described there seemed to be an infinite analogue of the Higman--Sims
graph which had played an important role in the researcher's thesis. The
encounter changed the course of the researcher's career. This paper
is the story of that event and its aftermath.

The final section of the paper suggests that Fra\"{\i}ss\'e classes of rigid
structures are a potentially interesting generalisation of Ramsey classes.
\end{abstract}

\section{Introduction}

Robert Woodrow was a good friend. He grew up in rural Canada, and I in rural
Australia, so we had a common background. We had some memorable adventures
together, and on my trip to Calgary to give the Louise and Richard Guy memorial
lecture, he looked after me very well, taking me to Banff and also to the
dinosaur-rich Alberta badlands. He is very much missed.

But, although we only have one joint paper~\cite{clptw}, he had a very
significant effect on my mathematical development. His paper~\cite{woodrow}
characterising Henson's homogeneous triangle-free graph was my introduction
to homogeneous structures and Fra\"{\i}ss\'e theory. So I am grateful to him
for launching me on this journey. In this paper I have mostly surveyed the
directions I have taken, with an emphasis on open questions.

For the remainder of this section, I provide some basic information about
homogeneous relational structures.

\subsection{Relational structures and homogeneity}

A relational structure is a structure over a first-order language with no
constant or function symbols. In other words, we have a number of ``named''
relations of various positive integer arities (the arity of a relation is the
number of arguments it takes). The commonest examples, each with a single
binary relation, are graphs and (totally or partially) ordered sets. We will
meet others later. The underlying set is called the \emph{domain} of the
structure. An \emph{isomorphism} between relational structures over the same
language is a bijection between the domains which preserves all the named
relations; an \emph{automorphism} of a structure is an isomorphism to itself.

If $M$ is a relational structure on a set $X$, and $Y$ is a subset of $X$, the
\emph{induced substructure} on $Y$ has the set $Y$, and relations which are the
restrictions to $Y$ of all the relations (in other words, we consider instances
where all the arguments belong to $Y$).

Note that we forbid relations with no arguments. Such a relation would be
either true or false in any given structure, and structures giving it different
truth values would not be comparable.

A structure $M$ is \emph{homogeneous} if every isomorphism between finite
induced substructures of $M$ can be extended to an automorphism of $M$. Note
that, in the early days, there was a confusion of terminology; some authors
(including Robert Woodrow) used the term
\emph{ultrahomogeneous} for this concept, and reserved \emph{homogeneous} for
the weaker concept that, if two finite induced substructures are isomorphic,
then some isomorphism between them can be extebded to an automorphism of $M$.
Note that the two notions here coincide for a special type of structure
which I will consider later, those in which every finite substructure is
\emph{rigid} (that is, has no non-trivial automorphisms). An important 
class of such structures consists of those for which one of the relations is
a total order, since there is a unique total order on a given finite set and
it is rigid. 

To try to minimise confusion, I will use the term ``weakly homogeneous'' for
the weaker concept, and if necessary use ``strongly homogeneous'' for the
stronger. I will discuss these further in the next section.

A variant of homogeneity was devised by Jaroslav Ne\v{s}et\v{r}il and me in
2004~\cite{cn}. A \emph{homomorphism} $f:A\to B$ between relational structures
is a map which maps tuples satisfying a relation in $A$ to tuples satisfying
the same relation in $B$. (Thus, in graphs, edges must map to edges; but
non-edges may map to edges, to non-edges, or to single vertices.) A structure
$A$ is said to have the \emph{HH property}, or to be
\emph{homomorphism-homogeneous}, if any homomorphism between finite induced
substructures of $A$ can be extended to an \emph{endomorphism} of $A$ (a
homomorphism from $A$ to itself). This concept has been studied and extended
by a number of authors~\cite{lt,ceg,appp}; but I will not discuss it further
here, apart from a brief mention in  Section~\ref{sec:R}.

\subsection{Graphs and groups}

My notation for graphs and groups is relatively standard. The \emph{complete
graph} $K_n$ with $n$ vertices has all possible edges between distinct vertices,
while the \emph{complete bipartite graph} $K_{m,n}$ has two sets with $m$ and
$n$ vertices, with all possible edges between the two sets. A \emph{complete
multipartite graph} is similar but with any number of sets.

The \emph{clique number} of a graph is the largest number of vertices in a
complete induced subgraph, and the \emph{independence number} is the largest
number of vertices in a set containing no edges. The \emph{line graph} of a
graph $\Gamma$ is the graph $L(\Gamma)$ whose vertices are the edges (the
pairs $\{x,y\}$ satisfying the adjacency relation), two vertices being joined
if the corresponding edges share a vertex.

Groups will often be automorphism groups of relational structures. I also need
some terms from permutation group theory. Permutations are written to the 
right of their arguments and composed left-to-right. Let $G$ be a permutation
group on a set $X$. The \emph{orbit} of an element $x\in X$ is the set of all
images of $x$ under elements of $G$. The group $G$ is \emph{transitive} if
there is a single orbit, that is, for any $x,y\in X$, there is an
element $g\in G$ with $xg=y$. It is \emph{primitive} if the only partitions
preserved by $G$ are the ``trivial'' ones, the partition into singletons and
the partition with a single part. For a positive integer $k$, the group is
$k$-transitive if, given any two $k$-tuples $(x_1,\ldots,x_k)$ and
$(y_1,\ldots,y_k)$ of distinct elements of $X$, there exists $g\in G$ with
$x_ig=y_i$ for $i=1,\ldots,k$. It is \emph{highly transitive} if it is 
$k$-transitive for all positive integers~$k$.

A permutation group is \emph{semiregular} if the stabiliser of any point is
the identity. (Topologists call such an action \emph{free}.) A group is
\emph{regular} if and only if it is semiregular and transitive. A regular
action of a group is isomorphic to the action on itself by right multiplication
(the \emph{right regular action}). There is also a \emph{left regular action},
where $g$ acts as left multiplication by $g^{-1}$; it is isomorphic to the
right regular action, and makes a brief appearance in Subsection~\ref{ss:sf}.

Two elements $g,h$ of a group $G$ are \emph{conjugate} if $h=x^{-1}gx$ for
some $x\in G$. This is an equivalence relation on $G$. Conjugate elements have
the same order, and the same cycle structure if $G$ is a permutation group.

The \emph{Classification of Finite Simple Groups} has had an enormous
effect on the theory of finite permutation groups. For just one instance, it
shows that a $6$-transitive finite permutation group must be a symmetric or
alternating group. (Things are different in the infinite case, as we shall
see.) I will refer to it briefly as CFSG.

\section{Finite homogeneous graphs}

The finite homogeneous graphs were classified by Gardiner~\cite{gardiner}:

\begin{theorem}
The finite homogeneous graphs are as follows:
\begin{enumerate}
\item the disjoint union of complete graphs of the same size;
\item the complements of the above (regular complete multipartite graphs);
\item the $5$-cycle;
\item the line graph of $K_{3,3}$.
\end{enumerate}
\end{theorem}

Here is a brief sketch proof. If $\Gamma$ is homogeneous and $v$ is any vertex,
then the induced subgraphs of $\Gamma$ on the neighbours and non-neighbours of
$v$ are both homogeneous. So we simply have to show that there is no further
homogeneous graph for which both of these subgraphs are on the list in the
theorem, and are the same for any $v$.

The theorem can be refined as follows. We say that a graph $\Gamma$ is
$t$-homogeneous (for a positive integer $t$) if any isomorphism between induced
subgraphs on at most $t$ vertices can be extended to an automorphism of
$\Gamma$. There is also a combinatorial version of this. We say that $\Gamma$
is \emph{$t$-tuple regular} if, given $t$-tuples $(x_1,\ldots,x_t)$ and
$(y_1,\ldots,y_t)$ of vertices such that the map $x_i\mapsto y_i$ for
$i=1,\ldots,t$ is an isomorphism of induced subgraphs, the number of common
neighbours of $x_1,\ldots,x_t$ is equal to the number of common neighbours of
$y_1,\ldots,y_t$. (This notion is called $t$-isoregularity in \cite{pech}.)
Thus $1$-tuple regular graphs are regular, while $2$-tuple
regular graphs are strongly regular.

The classification of $2$-homogeneous finite graphs follows from the
classification of rank~$3$ permutation groups, a consequence of CFSG
\cite{kl,liebeck,ls}.

With the algebraic methods that had been used to study the $3$-tuple regular
graphs~\cite{cgs}, I was able to prove~\cite{c:5tt}:

\begin{theorem}
A $5$-tuple regular finite graph is homogeneous.
\label{t:5tr}
\end{theorem}

There is more to the story. The \emph{Schl\"afli graph} on $27$ vertices is
given by Schl\"afli's ``double-six construction'' and describes the $27$ lines
in a general cubic surface~\cite{schlafli}. The \emph{McLaughlin graph} on
$275$ vertices was used by McLaughlin~\cite{mcl} to construct his sporadic
simple group.

For the Schl\"afli graph, there are three 4-vertex subgraphs which do not occur
as induced subgraphs: these are $K_4$ (one labelling), $K_3\cup K_1$ (four
labellings), and $K_4-e$ (six labellings). The remaining $53$ labelled graphs
on four vertices are single orbits of the automorphism group. So this graph
is $4$-homogeneous.

The McLaughlin graph is $4$-tuple regular, and contains every $4$-vertex graph
as an induced subgraph, but the independent sets of size $4$ fall into two
orbits under the automorphism group. So it is not $4$-homogeneous (but this is
the only failure of $4$-homogeneity).

From CFSG we know that there are no further $4$-homogeneous graphs. But the
existence of $4$-tuple but not $5$-tuple regular graphs beyond the Schl\"afli
and McLaughlin graphs is unknown.

\begin{qn}
Are there any further $4$-tuple regular graphs?
\end{qn}

A further open problem concerns recognition of these properties from the
\emph{spectrum} of a graph, the eigenvalues and multiplicities of the 
adjacency matrix. A graph property $P$ is \emph{determined by the spectrum}
if every graph which is cospectral with a graph having $P$ also has $P$.

\begin{theorem}
For $t\ne3$, the property of being $t$-tuple regular is determined by the
spectrum.
\end{theorem}

\begin{qn}
Is it true that $3$-tuple regularity is not determined by the spectrum?
\end{qn}

Now we can return to the question of weak and strong homogeneity, and note 
that they are the same for finite graphs. For a graph which is weakly
homogeneous is clearly $t$-tuple regular for all $t$, and so is homogeneous
by Theorem~\ref{t:5tr}. 

One can make a silly example where these concepts differ as follows. It is a
simple exercise to show that $(\mathbb{Q},{<})$ (the rational numbers with the
usual order) is homogeneous. Now define a relation $\rho$ on $\mathbb{Q}$ by
the rule that $\rho(x,y,z)$ holds if $x,y,z$ are all distinct and $x<y$. Then
the order on any set of size at least $3$ is determined by $\rho$; but there
are two $\rho$-isomorphisms between any two sets of size $2$, onlly one of
which extends to an automorphism.

\begin{qn}
Is there a weakly homogeneous infinite structure $M$ for which, for every $k$,
there is an isomorphism between two $k$-element substructures which does not
extend to an automorphism of $M$?
\end{qn}

For reasons of space I will say little about a more significant generalisation
of homogeneity for finite graphs with which I have recently been involved
\cite{bchk}. A graph $\Delta$ is an \emph{EPPA witness} for a graph $\Gamma$
if $\Gamma$ is an induced subgraph of $\Delta$ and every isomorphism between
induced subgraphs of $\Gamma$ can be extended to an automorphism of $\Delta$.
(Thus a graph is homogeneous if and only if it is its own EPPA witness.) The
name is an acronym for ``Extension Property for Partial Automorphisms''. 
Hrushovski~\cite{hrushovski} showed that every finite graph has a finite EPPA
witness, but interesting questions concerning the smallest EPPA witness remain.
The concept has important applications in showing significant properties of
infinite homogeneous graphs.

\section{Fra\"{\i}ss\'e's Theorem}

Clearly Gardiner's method cannot work for finding countable homogeneous
structures, since it is an induction on the cardinality of the structure.

In fact this question had been considered earlier by Roland
Fra\"{\i}ss\'e~\cite{fraisse}. It was Robert Woodrow's paper which introduced
me to Fra\"{\i}ss\'e's method, and I soon realised what a powerful construction
method it was.

The method depends on recognising a structure $M$ from its class of finite
substructures (more precisely, the finite structures over the same language
which are embeddable in $M$). Fra\"{\i}ss\'e calls this class the \emph{age}
of $M$. Clearly the age of a countable structure $M$ is a class $\mathcal{C}$
of finite relational structures and has the following properties:
\begin{enumerate}
\item It is closed under isomorphism.
\item It is \emph{hereditary}, that is, closed under taking induced 
substructures.
\item It contains (at most) countably many structures up to isomorphism.
\item It has the \emph{joint embedding property} (JEP): any two structures 
in $\mathcal{C}$ can both be embedded in a structure in $\mathcal{C}$.
\end{enumerate}

In future I will use the term ``class of relational structures'' to include
conditions (a), (b), (c) above. Note that (c) is automatic if the language
has only finitely many relation symbols.

It is not hard to see that the age of a countable homogeneous structure has
the \emph{amalgamation property} (AP): given structures
$A,B_1,B_2\in\mathcal{C}$ with embeddings $f_i:A\to B_i$ for $i=1,2$, there
exists $C\in\mathcal{C}$ and embeddings $g_i:B_i\to C$ such that
$f_1g_1=f_2g_2$, where $f_1g_1$ means $f_1$ followed by $g_1$.

Logicians normally forbid the empty set as a structure. If we allowed the
empty set, then (JEP) would be a consequence of (AP), since we have forbidden
relations of arity $0$.

Fra\"{\i}ss\'e's theorem~\cite{fraisse} asserts the following:

\begin{theorem}
\begin{enumerate}
\item
A class $\mathcal{C}$ of finite structures over a relational language $L$ 
(closed under isomorphism and induced substructures and at most countable up
to isomorphism) is the age of a countable homogeneous relational structure $M$
over $L$ if and only if it also satisfies (JEP) and (AP). 
\item
If $\mathcal{C}$ satisfies these conditions, then there is a unique countable
homogeneous structure $M$ (up to isomorphism) whose age is $\mathcal{C}$.
\end{enumerate}
\end{theorem}

A class $\mathcal{C}$ satisfying these conditions is a \emph{Fra\"{\i}ss\'e
class}, and the countable homogeneous structure $M$ is its \emph{Fra\"{\i}ss\'e
limit}.

For example, it is straightforward to see that the class of all finite graphs
is a Fra\"{\i}ss\'e class. Its Fra\"{\i}ss\'e limit is the subject of the
next section.

\medskip

I end this section with a question. A permutation group $G$ on a set $X$ is
\emph{oligomorphic} if it has only finitely many orbits on $X^n$ for all
positive integers $n$. By the theorem of Engeler, Ryll-Nardzewski and
Svenonius (Hodges~\cite[Theorem 6.3.1]{hodges}, a countable first-order
structure is $\aleph_0$-categorical (that is, specified uniquely up to
isomorphism among countable structures by its first-order theory) if and only
if its automorphism group is oligomorphic; and in this case, the automorphism
group orbits on $n$-tuples are precisely the $n$-types (where two $n$-tuples
have the same type if and only if they satisfy the same $n$-variable
first-order formulae). Further detail can be found in~\cite{hodges}.

A homogeneous structure over a finite relational language is oligomorphic,
since there are only finitely many isomorphism types of $n$-element
substructures for all $n$. The converse, however, is false.

Let $V$ be a vector space of countable dimension over the $2$-element field.
Let $X=V\setminus\{0\}$ and let $G$ be the group of all linear automorphisms
of $V$. Then $G$ is oligomorphic (for any finite set spans a finite subspace;
$G$ is transitive on subspaces of fixed finite dimension, and the stabiliser
of such a subspace has only finitely many orbits on its subsets). However,
there is no homogeneous relational structure over a finite language whose
automorphism group is $G$. This is because, if $k$ is the maximum arity of a
relation in the language, and $A$ and $B$ are subsets of cardinality $k+1$
such that $A$ is linearly dependent and $B$ satisfies the single linear
relation asserting that the sum of its elements is zero, then $A$ and $B$
induce isomorphic relational structures but are not in the same orbit of $G$.

\begin{qn}
Is there a property which distinguishes homogeneous structures over finite
languages among structures with oligomorphic automorphism groups?
\end{qn}

A necessary condition is given by the fact that the number of orbits on
$n$-sets of a homogeneous structure over a finite relational language is
bounded above by the exponential of a polynomial; there is no such bound for
$\aleph_0$-categorical structures. However the bound does hold in the above
example (the number of orbits on $n$-sets is roughly $2^{n^2/4}$).

\section{The Erd\H{o}s--R\'enyi random graph, or Rado's graph}
\label{sec:R}

In 1963, Erd\H{o}s and R\'enyi published a paper on the lack of symmetry of
finite random graphs. If a graph on $n$ vertices is formed by choosing edges
independently at random (that is, by tosses of a fair coin), then trivially
every $n$-vertex graph occurs with non-zero probability, and the probability
of a given graph is inversely proportional to the number of its automorphisms.
Slightly less trivially, the probability that the random graph has any 
non-identity automorphisms at all tends (rapidly) to $0$ as $n\to\infty$;
and further, if the distance between two graphs on the same vertex set is
measured by the number of insertions and deletions of edges required to change
one into the other, then with high probability a random graph lies at close to
the maximum possible distance from symmetry (that is, from a graph with a
non-identity automorphism).

By contrast, the paper contains a short tailpiece showing that a countably
infinite random graph has infinitely many automorphisms almost surely.

The real reason for this was given a decade later in the book \cite{es} by
Erd\H{o}s and Spencer on the probabilistic method. There is a particular
countable graph $R$ which has the property that a random countable graph is
isomorphic to $R$ almost surely. Moreover, $R$ has infinitely many
automorphisms; indeed, it is homogeneous. 

The proof is remarkably simple. Consider the ``extension property'' that,
given two disjoint finite sets $U$ and $V$ of vertices, there is a vertex $z$
joined to every vertex in $U$ and to none in $V$. Easy arguments show that
the countable random graph has this property almost surely. The back-and-forth
method from model theory shows that two countable graphs with this property are
isomorphic. Moreover, it is clear that the age of a graph with this property
consists of all finite graphs. Also, given an isomorphism between finite
subgraphs, take two enumerations of the graph in which these finite subgraphs
occur first, and then use back-and-forth to extend the isomorphism between them
to an automorphism.

I will denote this graph by $R$, and call it the (countable) random graph.
Erd\H{o}s and Spencer say that this result ``demolishes the theory of
countable random graphs''; I would say that, in fact, it creates the new theory
of the graph $R$.

Erd\H{o}s and Spencer gave no explicit construction of $R$, since their
argument is a textbook example of a non-constructive existence proof (something
which occurs with probability~$1$ certainly exists). But in the intervening
decade, Rado~\cite{rado} had given an explicit construction. The vertex set is
the set $\mathbb{N}$ of natural numbers. Given two natural numbers $x$ and $y$,
we may assume that $x<y$; then join $x$ and $y$ if the $x$th digit in the
base-$2$ expression for $y$ is $1$. It is a simple exercise to verify the
``extension property'' for this graph. 

Note that, if instead we take this
relation to be directed, from smaller to larger, we obtain a model of
\emph{hereditarily finite set theory}, that is, Zermelo--Fraenkel set theory
with the Axiom of Infinity replaced by its negation. More generally, if we take
a countable model of the ZFC axioms Empty Set, Pairing, Union and Foundation,
and ignore directions, the undirected graph we obtain is $R$.

\begin{qn}
At what point between 1963 and 1974 did Erd\H{o}s realise the result which now
appears in \cite{es}? Did Rado know about the result of Erd\H{o}s and R\'enyi,
or vice versa?
\end{qn}

Another explicit construction of $R$ is the following. The vertex set is the
set $\mathbb{P}_1$ of primes congruent to $1$~(mod~$4$); join $p$ to $q$ if
$q$ is a quadratic residue mod~$p$. By quadratic reciprocity, the graph is
undirected; using the Chinese Remainder Theorem and Dirichlet's theorem on
primes in arithmetic progressions, we can verify the extension property which
characterises $R$.

The graph $R$ was a rich field for study. One of the first questions I asked
myself was, does it have an automorphism permuting all the vertices in a
single cycle? If there is such an automorphism $\sigma$, we can label the
vertices with the integers so that $\sigma$ acts as the cyclic shift.
Let $S$ be the set of positive neighbours of the vertex $0$. Then the whole
graph is determined by $S$: we have $x$ joined to $y$ if and only if
$|y-x|\in S$. I was able to show that, if we choose $S$ at random by tossing
a fair coin, then the graph we obtain is isomorphic to $R$ almost surely.
Moreover, two cyclic automorphisms of $R$ are conjugate in the full
automorphism group if and only if the corresponding sets $S$ are equal. This
demonstrates that $\Aut(R)$ has $2^{\aleph_0}$ conjugacy classes of cyclic
automorphisms; so in particular, its cardinality is $2^{\aleph_0}$.

There is another approach to this. The countable graph is specified by a
sequence of $0$s and $1$s (for non-edges and edges); we can regard the set of
all such sequences as a probability space. But we can also regard it as a
metric space, where the distance between two sequences is $1/2^n$ if they first
disagree in the $n$th position. The metric space is complete, and so the notion
of Baire category gives a second interpretation of ``almost all'': the large
sets are \emph{residual}, that is, they contain a countable intersection of
dense open sets. (The Baire category theorem asserts that residual sets are
non-empty, and indeed they meet any open set non-trivially.)

The set of graphs isomorphic to $R$ is residual in the space, again giving a
nonconstructive existence proof. On the face of it, probability is more
informative, since sets may have probability between $0$ and $1$. But Baire
category provides a tool which is much easier to use. In the space of all
binary sequences, a set $S$ is open if and only if it is \emph{finitely
determined} (that is, any member of $S$ has an initial finite sequence all of
whose continuations lie in $S$) and dense if and only if it is \emph{always
reachable} (any finite sequence is an initial subsequence of a member of $S$).
A residual set contains an intersection of countably many sets with these
properties.

Now it is easy to see that the ``extension property'' above for fixed $U$ and
$V$ is finitely determined and always reachable, and there are only countably
many choices for $U$ and $V$; so graphs isomorphic to $R$ form a residual set
of all graphs on a countable vertex set.

I conclude with two properties of $R$ from \cite{cn}. A graph $\Gamma$ is a
\emph{spanning subgraph} of a graph $\Delta$ if they have the same vertex set
and every edge of $\Gamma$ is an edge of $\Delta$.

\begin{theorem}
\begin{enumerate}
\item A countable graph $\Gamma$ contains $R$ as a spanning subgraph if and
only if any finite set of vertices of $\Gamma$ has a common neighbour.
\item A graph containing $R$ as a spanning subgraph is 
homomorphism-homogeneous.
\end{enumerate}
\end{theorem}

\section{Countable homogeneous graphs}

Robert Woodrow's paper~\cite{woodrow} shows that there are precisely four
countable homogeneous graphs without triangles: a countable null graph; a
countable disjoint union of edges; a complete bipartite graph with two
countably infinite parts; and one further example containing all finite
triangle-free graphs. This graph can be explained as the Fra\"{\i}ss\'e limit
of the class of all finite triangle-free graphs.

In fact, this last graph was the first of an infinite family of countable
homogeneous graphs. For $k\ge3$, \emph{Henson's graph} $H_k$ is the
Fra\"{\i}ss\'e limit of the class of all graphs containing no induced subgraph
$K_k$ (complete on $k$ vertices). The proof of the amalgamation property is
simple: make the amalgam of $B_1$ and $B_2$ over $A$ by making no
identifications other than the vertices in $A$ and putting no extra edges
between $B_1$ and $B_2$. (Henson's construction~\cite{henson} was different:
he built his graphs inductively inside $R$.)

A little later, Woodrow with his doctoral adviser Alistair Lachlan gave a
complete list~\cite{lw}:

\begin{theorem}
A countable homogeneous graph is one of the following:
\begin{enumerate}
\item A disjoint union of complete graphs of the same size, where either the
number of parts or the size of the parts (or both) is countably infinite.
\item The complement of a graph in (a).
\item A Henson graph $H_k$ for $k\ge3$.
\item The complement of a graph in (c).
\item The graph $R$.
\end{enumerate}
\end{theorem}

Several further classifications of homogeneous relational structures have been
found, but this stands as probably the most significant.

As I have said, Woodrow's characterisation of $H_3$ particularly appealed to
me. In my doctoral thesis, the graph used by Higman and Sims to construct a new
sporadic simple group~\cite{hs} plays a starring role. At the time, we called
it the Higman--Sims graph; but it was pointed out later that it had been
constructed some time earlier by the statistican Dale Mesner~\cite{mesner1},
and later proved unique by him~\cite{mesner2} (though he did not think to
examine its automorphism group).

This graph has the properties that its automorphism group acts primitively on
vertices, and is transitive on ordered edges and ordered non-edges; the graph
has no triangles, and the stabiliser of a vertex is $3$-transitive on the set
of neighbours of that vertex. For comparison, the automorphism group of $H_3$
is primitive on vertices and transitive on ordered edges and ordered non-edges;
the stabiliser of a vertex is \emph{highly transitive} (that is,
$k$-transitive for every positive integer $k$) on the set of neighbours of that
vertex. This highly transitive group was like none I had seen before, and I
was inspired to look for further examples.

\section{B-groups}

Wielandt~\cite[Definition 25.1]{wielandt} defined a \emph{B-group} to be a
group with the following property: any primitive permutation group which
contains the group $G$ in its right regular action must be $2$-transitive.
(In other words, if we add permutations to kill all the non-trivial
$G$-invariant equivalence relations, then we kill all non-trivial $G$-invariant
binary relations.)

The groups were named for William Burnside, who proved that a finite cyclic
group of order a proper power of a prime number has this property. (We cannot
call these groups \emph{Burnside groups}, since this term was already
established for finitely generated infinite groups of finite exponent.) 

The study of these groups was very influential, and led to important 
developments in the theory of Schur rings. However, the Classification of
Finite Simple Groups completely changed the picture. For example, we know that
for almost all natural numbers $n$ (all but $O(x/\log x)$ numbers
below~$x$), the only primitive groups of degree $n$ are the symmetric and
alternating groups, and so every group of order $n$ is a B-group~\cite{cnt}.

However, the picture in the infinite case is very different. Ken Johnson visited
Oxford in 1981--2 and asked whether there are any countably infinite B-groups.
Graham Higman was able to give a rather general condition ensuring that a
countable group was not a B-group. Then Johnson and I~\cite{cj} showed that
every group satisfying a weaker form of Higman's condition is a regular
subgroup of the automorphism group of $R$. Since $\Aut(R)$ is primitive but
not $2$-transitive, all groups satisying this condition fail to be B-groups.

Here is a statement of our theorem. In a group $G$, a \emph{square-root set}
is a set of the form $\sqrt{a}=\{x:x^2=a\}$ for a fixed element $a$; it is
\emph{non-principal} if $a\ne 1$. 

\begin{theorem}
Let $G$ be a countable group which cannot be written as the union of finitely
many translates of non-principal square-root sets together with a finite set.
Then $G$ is embedded as a regular subgroup of $\Aut(R)$.
\end{theorem}

For any finite or countable group $G$, the group $G\times C_\infty$ satisfies
the hypotheses of this theorem. Thus $\Aut(R)$ embeds every finite or 
countable group as a semiregular subgroup.

As far as I know, no countable group has ever been shown to be a B-group. The
best candidate I know is the infinite dicyclic group
\[\langle x,y\mid y^4=1,y^{-1}xy=x^{-1}\rangle\]
which fails to satisfy our version of Higman's condition: it can be written as
$\sqrt{y^2}\cup y\sqrt{y^2}$.

\begin{qn}
Is there a countable B-group?
\end{qn}

\section{Homogeneous Cayley objects}

Although I had to leave the question of countable B-groups unresolved, further
work on it let me to interesting mathematics and collaborations with two
remarkable mathematicians.

The use of $R$ in the preceding section suggested turning the question on its
head. Let $M$ be a countable structure. I will say that $M$ is a \emph{Cayley
object} for the countable group $G$ if its point set can be identified with $G$
so that the right regular action of $G$ is
contained in the automorphism group of $M$. This is motivated by the notion of
\emph{Cayley graph}, a graph which is a Cayley object for a group $G$.
Introduced by Cayley in the 19th century, these graphs are now fundamental
objects in geometric group theory. My hope was that, by choosing homogeneous
structures whose automorphism group was primitive but not $2$-transitive and
expressing them as Cayley objects for suitable groups, I would find further
non-B-groups. Although I had no success with this, the search did lead me
into some very interesting byways, some of which I will now describe.

\subsection{Hypergraphs}

For $k>2$, there is a unique countable homogeneous universal $k$-uniform
hypergraph $H$. It can be shown that it is a Cayley object for every countable
group. However, this gives no further B-groups, since the automorphism group
of $H$ is $(k-1)$-transitive. (Thus all degrees of transitivity are realised
by infinite permutation groups.)

\subsection{Sum-free sets}
\label{ss:sf}

Of course, one of the first examples I looked at was $H_3$, Henson's
triangle-free graph. I looked first at cyclic automorphisms.
Henson~\cite{henson} had shown in his original paper that the graph does admit
cyclic automorphisms. Following the approach that worked for $R$, what was
needed was to characterise the sets $S$ of positive integers such that the
graph $\Gamma(S)$ obtained by joining two integers $x$ and $y$ if $|y-x|\in S$
is isomorphic to Henson's graph. It is easy to see that $\Gamma(S)$ is
triangle-free if and only if $S$ is \emph{sum-free}, that is, does not
contain $x,y,z$ with $x+y=z$. (We allow $x=y$ here.) So the question was: 
do almost all sum-free sets $S$ produce Henson's graph? It is not hard to write
down a condition on $S$ for $\Gamma(S)$ to be Henson's graph; such sum-free
sets are called \emph{universal}.

I tried and failed to find a nice isomorphism-invariant measure to do the job.
Later, by completely different techniques resembling graphon theory, Petrov
and Vershik~\cite{pv} found such a measure, and this was extended to a much
wider class of structures by Ackerman, Freer and Patel~\cite{afp}.

However, Baire category was able to do the job. It is not hard to show that
the set of universal sum-free sets is residual in the set of all sum-free sets.
The explicit construction of such a set involved leaving increasingly long
gaps, and I conjectured that a universal sum-free set has density zero. This
conjecture was confirmed by Schoen~\cite{schoen}.

There is a nice pattern here. Two old results related to Ramsey's theorem are:
\begin{itemize}
\item Van der Waerden's theorem~\cite{vdw}: in a partition of $\mathbb{N}$
into finitely many parts, one part contains arbitrarily long arithmetic
progressions.
\item Schur's theorem~\cite{schur}: in a partition of $\mathbb{N}$ into
finitely many parts, one part contains two numbers and their sum.
\end{itemize}
Famously, van der Waerden's theorem was strengthened by Szemer\'edi~\cite{szem}
to the statement that a set of natural numbers of positive upper density
contains arbitrarily long arithmetic progressions. But the analogous
strengthening of Schur's theorem is false: the set of odd numbers is sum-free
but has density $\frac{1}{2}$. The results noted above give a kind of
replacement: almost all sum-free sets in the sense of Baire category (the
universal ones) are excluded by the positive upper density condition.

Sum-free sets proved to be a fascinating topic. Of course, I wanted to count
them, and this problem which I posed in 1987~\cite{c:sf} caught the interest of 
Paul Erd\H{o}s, and led to my first joint paper with him~\cite{ce}.
The conjecture was that the number of sum-free subsets of $\{1,\ldots,n\}$ is
asymptotically $c_e2^{n/2}$ if $n$ is even and $c_o2^{n/2}$ if $n$ is odd,
where $c_e$ and $c_o$ are two constants (whose values are roughly $6.8$ and
$6.0$). More precisely, the conjecture asserts that almost
all such sets either consist of odd numbers or have all entries in 
$(n/2-w(n),n]$, where $w(n)$ is a function growing arbitrarily slowly. We were
able to show that for these two types our asymptotic was correct. Later
Green~\cite{green} and Sapozhenko~\cite{sapozhenko} proved the conjecture, by
showing that the number of sets not of these two types was asymptotically
smaller than $2^{n/2}$.

\begin{qn}
Estimate the constants $c_e$ and $c_o$ more precisely. Are they transcendental?
\end{qn}

There is a probability measure defined on sum-free sets of natural numbers as
follows. Start with the empty set, and consider the natural numbers in turn.
If $n=x+y$ where $x,y\in S$, then $n\notin S$; otherwise include $n$ in $S$
by the toss of a fair coin. This measure has some interesting properties but
little is known about it. Here is a sample question, suggested by experiment.

\begin{qn}
Examine this space further. Is it true that a random sum-free set constructed
as above has a density almost surely, and that the spectrum of densities is
discrete above $\frac{1}{6}$ but has a continuous part below this value?
\end{qn}

It is known, for example, that the probability that a random sum-free set 
contains no even numbers is non-zero (it is approximately $0.218$), and
that conditioned on this, the density is almost surely $\frac{1}{4}$.
Other types of sum-free set (for example, those contained in the residue
classes $1$ and $4$ (mod~$5$), or those in which $2$ is the only even number)
also have positive probability~\cite{c:rsf,cc}.
\medskip

Henson's graph $H_3$ is a Cayley graph for various countable groups, but none
which are not already accounted for by the random graph $R$.

What about the higher Henson graphs? Henson had shown that, for $k\ge4$, the
graph $H_k$ has no cyclic automorphism. I extended this to show that it is
not a \emph{normal Cayley graph} for any countable group. (A normal Cayley
graph is a graph $\Cay(G,S)$ for which $S$ is a normal subset of $G$, that is,
fixed by conjugation; equivalently, the graph admits both the right and the
left regular actions of $G$.) In particular, $H_k$ is not a Cayley graph for
any abelian group if $k\ge4$.

Then Cherlin~\cite{cherlin} showed that these graphs are Cayley graphs for
the free group of countable rank, and also for the free
nilpotent group of class~$2$ and countable rank. His techniques also
work for a wider class of graphs, including some related to the Urysohn space,
as well as for various directed graphs and metric spaces.

\subsection{The Urysohn space}

After I talked about $R$ at the European Congress of Mathematics in Barcelona
in 2000, Anatoly Vershik introduced himself to me and told me about the
\emph{Urysohn space}, the universal homogeneous Polish space (complete
separable metric space). It is too big to be the Fra\"{\i}ss\'e limit of finite
metric spaces, but as Urysohn realised in the 1920s, there is a simple way
around this. A \emph{rational metric space} is a metric space in which all the
distances are rational numbers. There are only countably many rational metric
spaces, and they do form a Fra\"{\i}ss\'e class; their Fra\"{\i}ss\'e limit is
the \emph{rational Urysohn space}, and the usual process of completion now
gives the Urysohn space.

With Vershik's guidance, and the experience of dealing with other
Fra\"{\i}ss\'e classes, we were able to prove some striking results about
isometries of Urysohn space, especially abelian groups acting regularly on it.
(This is done by finding a cyclic automorphism of the rational Urysohn space;
its orbits on the completion are dense, so its closure is a transitive abelian
group.  These give a number of ways of putting an abelian group structure on
the Urysohn space (although the constructions are not explicit).
See~\cite{cv}.

\begin{qn}
Noting that the completion process may introduce torsion, we can ask: can we
say anything about the torsion subgroups of the abelian groups which arise?
In particular, is there a torsion-free transitive abelian subgroup?
\end{qn}

\subsection{Multiorders}

A little thought shows that the structure $(\mathbb{Q},{<})$ (the rational
numbers as ordered set) is homogeneous (its age is the set of all finite
total orders), and is a Cayley object for the additive group of $\mathbb{Q}$.
Indeed, this is the example which Fra\"{\i}ss\'e generalised in his
paper~\cite{fraisse}.

A \emph{multiorder} is simply a set carrying a number of total orders. If there
are $m$ orders, we speak of an $m$-order.

A finite $2$-order on $n$ points is simply a permutation of $\{1,\ldots,n\}$,
regarded as a reordering of these labels, not as a function on the set of 
labels. The first order can be used to label the $n$ points with $1,2,\ldots,n$;
then the second order rearranges the labels. Now a substructure of a 
$2$-order (as relational structure) gives a permutation on $\{1,\ldots,k\}$ in
which the labels come in the same order as in the permutation on
$\{1,\ldots,n\}$. For example, $[1,3,2]$ is a substructure of $[2,4,1,3]$,
using the first, second and fourth entries. This is exactly the notion of
subpermutation which occurs in the theory of permutation patterns~\cite{wilf}.

The finite $m$-orders form a Fra\"{\i}ss\'e class. I investigated its
Fra\"{\i}ss\'e limit as a Cayley object~\cite{c:ajc}. All I was able to prove
was that the homogeneous $m$-order is a Cayley object for the free abelian group
$(C_\infty)^n$ if and only if $m<n$. The proof required Kronecker's theorem
on Diophantine approximation~\cite{kronecker}.

For example, take the elements of $(C_\infty)^2$ as points of the plane with
integer coordinates. Now to construct a dense total order on this set, take
a line with irrational slope, and slide it across the plane. The order of the
points is given by the order in which they are hit by the sliding line.
Translations of the set by vectors with integer coordinates preserve the order.

A group is a Cayley object for a total order if and only if it is \emph{right
ordered}; for the order to be homogeneous we require the right order on the
group to be dense.

\begin{qn}
Find other examples of groups for which the universal $m$-order is a Cayley
object, for $m>1$.
\end{qn}

\subsection{Miscellanea}

There are interesting graphs which are in a sense ``almost homogeneous'': we
can produce a homogeneous structure by adding a relation. One example is the
generic \emph{bipartite graph}: this is not homogeneous, since non-adjacent
vertices may be at distance $2$ or $3$, but if we add an equivalence relation
to the language and interpret it as the bipartition, we obtain a homogeneous
structure. But its automorphism group is imprimitive, so we obtain no B-groups.

A more interesting example is Covington's N-free graph~\cite{covington}. A
graph is \emph{N-free}, or a \emph{cograph}, if it contains no path of
length~$4$ as induced subgraph. If we have a $3$-clique $C$ in such a graph,
it is easy to see that if two vertices $v,w$ are each joined to a single
vertex in $C$, then it must be the same vertex. Thus, we need a relation to
pick out one of the three vertices as a potential single neighbour of a point
outside. A similar remark holds for independent sets of size $3$. Any other
$3$-vertex subgraph already has a distinguished vertex. Covington showed that
by adding a ternary relation which distinguishes one of its three
arguments suitably, we obtain a homogeneous structure. This structure is a
Cayley object for the countable elementary abelian $2$-group.

This ternary relation is called a \emph{C-relation}; it plays a role in the
classification of Jordan groups of countable degree~\cite{an}, and will crop
up again in the last section of this paper, where I will give a description of
it.

\begin{qn}
Is it true that a group for which Covington's graph is a Cayley object must be
a $2$-group? Are there examples other than the elementary abelian group?
\end{qn}

\section{Reducts and overgroups}

Let $M$ be a homogeneous relational structure. A structure $N$ on the same
point set is a \emph{reduct} of $M$ if the relations of $N$ have first-order
definitions without parameters in $M$. Reducts are defined up to equivalence, 
two structures equivalent if each is a reduct of the other. As an example,
here are the reducts of the structure $(\mathbb{Q},{<})$. This list is a
consequence of my first paper on infinite structures~\cite{c:trans}, though I
did not know at the time that I was proving that.
\begin{enumerate}
\item The structure $(\mathbb{Q},{<})$.
\item The \emph{betweenness relation} on $\mathbb{Q}$, the ternary relation
$\beta$, where $\beta(x,y,z)$ if and only if $x<y<z$ or $z<y<x$.
\item The \emph{circular order} $\gamma$, the ternary relation which holds
if and only if $x<y<z$ or $y<z<x$ or $z<x<y$.
\item The \emph{separation relation} $\sigma$, the quaternary relation for
which $\sigma(w,x,y,z)$ holds if and only if $w$ and $y$ are between $x$
and $z$ and \emph{vice versa}.
\item The pure set $\mathbb{Q}$ with no relations.
\end{enumerate}

There is a natural topology on the symmetric group, the topology of pointwise
convergence; a basis for the open neighbourhoods of the identity is given by
the pointwise stabilisers of finite tuples. In this topology, a subgroup of
the symmetric group is closed if and only if it is the automorphism group of
a relational structure, which can be taken to be homogeneous (but we cannot
assume that the language is finite). Now it is clear that, if $N$ is a reduct
of $M$, then $\Aut(N)$ is a closed overgroup of $\Aut(M)$, and conversely
(for homogeneous structures over a finite language).

I suggested a list of reducts of $R$, and Simon Thomas~\cite{thomas} proved
that it was complete. I will describe the groups rather than the structures.
\begin{enumerate}
\item $\Aut(R)$.
\item The group of automorphisms and anti-automorphisms (isomorphisms to the
complement) of $R$.
\item The group of switching automorphisms of $R$. (Switching is the operation
of a graph given by choosing a subset $Y$ of the vertex set $X$, and
exchanging edges and non-edges between $Y$ and its complement, keeping the
induced subgraphs on $Y$ and on $X\setminus Y$ as they were; a switching
automorphism is a bijection on the vertex set which maps the graph to an
image under switching.)
\item The group of switching automorphisms and anti-automorphisms of $R$.
\item The symmetric group on the vertex set of $R$.
\end{enumerate}

Note the structural similarity with the reducts of $(\mathbb{Q},{<})$. In 
particular, the degree of transitivity of the reducts (a)--(d)
is $1$, $2$, $2$, $3$, and (a) and (c) are normal subgroups of (b) and (d) of
index $2$, in each case. This result led Thomas to conjecture
that a homogeneous structure over a finite relational language has only
finitely many reducts. This has been verified in a number of cases.

\begin{qn}
Decide Thomas' conjecture.
\end{qn}

The automorphism group of a homogeneous relational structure may have many
overgroups in the symmetric group which are not closed, and hence not
reducts. For example, the finitary symmetric group $F$ on a countable set $X$
(the group of permutations which move only finitely many points)
is a normal subgroup of the symmetric group, so for any subgroup $G$, the
product $FG$ is a subgroup; it is highly transitive, so its closure is the
symmetric group.

But some of these overgroups are interesting in their own right. My only 
joint paper with Robert Woodrow came about when he and his colleagues
Claude Laflamme and Maurice Pouzet were examining automorphism groups of
hypergraphs (with infinite hyperedges, so not relational structures: for
example, the hypergraph $\mathcal{H}$ whose edges are the subsets inducing
copies of $R$), while Sam Tarzi and I were looking at homeomorphism groups of
filters or topologies (for example, the neighbourhood filter of $R$, generated
by the vertex  neighbourhoods in $R$), and groups of permutations which are
``almost automorphisms'' of $R$. We combined forces to produce the paper
\cite{clptw}. The paper contains a couple of open questions about the
relationships between these groups.

Examples of groups of ``almost automorphisms'' include $\Aut_1(R)$, the group
of permutations changing only finitely many adjacencies (edges to non-edges or
\emph{vice versa}); $\Aut_2(R)$, the group of permutations changing only 
finitely many adjacencies at each vertex; and $\Aut_3(R)$, the group of
permutations changing only finitely many adjacencies at all but possibly 
finitely many vertices.

The paper contains a couple of open questions about the relationships between
these groups. For example, let $\mathrm{FAut}(\mathcal{H})$ denote the group
of permutations $g$ such that there is a finite subset $S$ of $R$ such that,
for every edge $E$ of $\mathcal{H}$, both $(E\setminus S)g$ and
$(E\setminus S)g^{-1}$ are edges of $\mathcal{H}$. It is known that 
\[\Aut(\mathcal{H})\cdot\mathrm{FSym}(R)\le\mathrm{FAut}(\mathcal{H}),\]
where $\mathrm{FSym}(R)$ is the group of permutations of finite support on the
vertex set of $R$, but it is not known whether equality holds.

\begin{qn}
Resolve the outstanding questions about equalities and inclusions in 
\cite{clptw}.
\end{qn}

The paper also shows:

\begin{theorem}
An overgroup of $\Aut(R)$ which is not a reduct must be either contained in
the group $B$ of switching automorphisms and anti-automorphisms, or highly
transitive.
\end{theorem}

For, if it is not contained in $B$, then its closure must be the symmetric
group.

\begin{qn}
Determine the overgroups of $\Aut(R)$ contained in $B$.
\end{qn}

An example is the group of switching automorphisms where the switching set is
finite.

\section{Ramsey classes}

A simple form of Ramsey's theorem states that, given positive integers $a$ and
$b$, there exists $c$ such that, if the $a$-element subsets of a $c$-element
set are coloured red and blue, there will be a monochromatic $b$-element set.

In the 1980s, interest grew in a structural version of the theorem, where we
replace sets by structures in a class $\mathcal{C}$ over a finite relational
language. As usual, we assume our classes are closed under isomorphism and
hereditary. We  also have to replace substructures by embeddings, since
different embeddings with the same image destroy the Ramsey property.

Thus, $\mathcal{C}$ is a \emph{Ramsey class} if, given structures
$A,B\in\mathcal{C}$, there exists $C\in\mathcal{C}$ such that, if the 
embeddings $A\to C$ are coloured red and blue, there exists an embedding
$B\to C$ such that all the embeddings $A\to C$ with image in $B$ have the same
colour.

Jaroslav Ne\v{s}et\v{r}il~\cite{nesetril} showed two general properties of
Ramsey classes.

\begin{theorem}
Let $\mathcal{C}$ be a Ramsey class over a finite relational language. Then
\begin{enumerate}
\item $\mathcal{C}$ is a Fra\"{\i}ss\'e class;
\item if the structures in $\mathcal{C}$ are non-trivial (not pure sets), then 
they are rigid (that is, they have trivial automorphism groups).
\end{enumerate}
\end{theorem}

Ne\v{s}et\v{r}il found many examples of Ramsey classes, and proposed a programme
to classify them: determine the Fra\"{\i}ss\'e classes, and then decide which
of them have the Ramsey property.

All of Ne\v{s}et\v{r}il's examples were rigid because a total order was part of
the structure: more formally, they have a total order as a reduct (which may
be done by imposing a total order on all the structures in the family). At the
time, I constructed a Fra\"{\i}ss\'e class of rigid structures by a different
method, and wondered whether or not it was a Ramsey class. 

The answer came in the 1990s with the theorem of Kechris, Pestov, and
Todor\v{c}evi\'c~\cite{kpt}. It uses the topology of pointwise convergence on
the symmetric group, introduced in the preceding section. A topological
group $G$ is said to be \emph{extremely amenable} if, for any continuous
action of $G$ on a compact space $X$, there is a point of $X$ fixed by $G$.

\begin{theorem}
Let $\mathcal{C}$ be a Fra\"{\i}ss\'e class, with Fra\"{\i}ss\'e limit $M$.
Then $\mathcal{C}$ is a Ramsey class if and only if $\Aut(M)$ is extremely
amenable.
\end{theorem}

This answers my question. The set of total orders on a countable set can be
shown to be compact, in a natural  topology, with the symmetric group acting
continuously on it. So, if $\mathcal{C}$ is a Fra\"{\i}ss\'e class with
Fra\"{\i}ss\'e limit $M$, then $\Aut(M)$ fixes a total order. It follows from
the Engeler--Ryll-Nardzewski--Svenonius theorem of model theory that this total
order is a reduct of $M$, and so can be defined without parameters. Hence every
finite structure in $M$ carries a total order, and thus is rigid.

I will describe my example in the next section, where I will show that it does
not have a total order as a reduct, and hence is not a Ramsey class. With
Siavash Lashkarighouchani~\cite{cl}, I was able to find an explicit failure of
the Ramsey property in my class, and indeed in any Fra\"{\i}ss\'e class of
rigid structures which is not a Ramsey class. This failure involves structures
$A$ and $B$ with $|A|=2$.

\section{Fra\"{\i}ss\'e classes of rigid structures}

This final section considers Fra\"{\i}ss\'e classes of rigid structures, which
(as noted in the Introduction) are precisely those for which the weak and
strong versions of homogeneity (formerly called homogeneity and
ultrahomogeneity) coincide. This is a class which might repay investigation,
and I suggest some pointers.

A Fra\"{\i}ss\'e class $\mathcal{C}$ is said to have the \emph{strong
amalgamation property} if, for any amalgamation $A\to B_1,B_2\to C$, the
structure $C$ and the embeddings of $B_1$ and $B_2$ into $C$ can be chosen
so that the intersection of the images of $B_1$
and $B_2$ is precisely the image of $A$ (and not larger). This is equivalent
to saying that, if $M$ is the Fra\"{\i}ss\'e limit and $G=\Aut(M)$, the 
stabiliser of a finite set $F$ in $G$ has no finite orbits outside $F$. In
model-theoretic terms, this says that algebraic closure in $M$ is trivial.

Classes with the strong amalgamation property are important in several areas.
For example, the theorem of Ackerman \emph{et al.}~\cite{afp} on invariant
measures shows that this condition is necessary and sufficient for their
construction of an invariant measure. Also, if $G$ is a permutation group on
a countable set
in which the stabiliser of a finite set has no finite orbits outside the set,
then either $G$ is highly transitive, or $G$ preserves a non-trivial topology
(and if primitive, it preserves a non-trivial filter)~\cite[Section 4]{c:top}.

If two Fra\"{\i}ss\'e classes $\mathcal{C}_1$ and $\mathcal{C}_2$ have the
strong amalgamation property, we may form the class
$\mathcal{C}=\mathcal{C}_1\odot\mathcal{C}_2$ over the union of the languages
for the two classes: a structure in this class consists of a finite set with
independently-chosen structures from $\mathcal{C}_1$ and $\mathcal{C}_2$. It
is straightforward to see that $\mathcal{C}$ is also a Fra\"{\i}ss\'e class
with strong amalgamation.

For example, the class $\mathcal{C}_m$ of $m$-orders can be defined inductively
by the rules that $\mathcal{C}_1$ is the class of total orders and
$\mathcal{C}_m=\mathcal{C}_{m-1}\odot\mathcal{C}_1$ for $m>1$.

A \emph{tournament} is a set with a binary relation $\tau$ such that
$\tau(x,x)$ never holds, while for $x\ne y$, exactly one of $\tau(x,y)$ and
$\tau(y,x)$ holds. It is easy to see that the class of finite tournaments is
a Fra\"{\i}ss\'e class with strong amalgamation. Its Fra\"{\i}ss\'e limit is
the \emph{generic tournament}, with properties resembling those of the random
graph; in particular, it can be constructed using the set $\mathbb{P}_{-1}$
of primes congruent to $-1$~(mod~$4$), with an arc from $p$ to $q$ if $q$ is
a quadratic residue mod~$p$.

Moreover, if $T$ is a finite tournament, then $|\Aut(T)|$ is odd: for, if it
were even, then it would contain an element of order~$2$, which would
interchange two vertices; but this is clearly not possible.

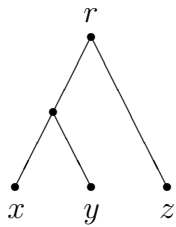
\begin{figure}[htbp]
\begin{center}
\setlength{\unitlength}{1mm}
\begin{picture}(20,30)
\multiput(10,25)(-5,-10){3}{\circle*{1}}
\multiput(10,5)(10,0){2}{\circle*{1}}
\put(10,25){\line(1,-2){10}}
\put(10,25){\line(-1,-2){10}}
\put(5,15){\line(1,-2){5}}
\put(9,27){$r$}
\put(-1,1){$x$}
\put(9,1){$y$}
\put(19,1){$z$}
\end{picture}
\caption{\label{f:bt}C-relation $\gamma(x,y;z)$}
\end{center}
\end{figure}

The other ingredient is a C-relation, one of those involved in the 
classification of infinite Jordan groups~\cite{an}. Let $X$ be the set of
leaves of a rooted binary tree. Define the ternary relation $\gamma$ on $X$
by the rule that $\gamma(x,y;z)$ holds if
\[[x,r]\cap[y,r]\supset[x,r]\cap[z,r]=[y,r]\cap[z,r],\]
where $r$ is the root, and $[a,b]$ denotes the unique path from $a$ to $b$.
See Figure~\ref{f:bt}.

It can be shown that the C-relation on the leaves determines the binary tree,
and hence that the class of finite C-relations is a Fra\"{\i}ss\'e class with
the strong amalgamation property. Moreover, the automorphism group of any
C-relation is a $2$-group. For such an automorphism must preserve the binary
tree. An automorphism of odd prime order fixes the root, and hence its two 
children, and so on down the tree; thus no such automorphism can exist. So
the automorphism group of any finite C-structure is a $2$-group.

Now, if $\mathcal{C}_1$ is the class of finite tournaments, and $\mathcal{C}_2$
the class of finite C-relations, then
$\mathcal{C}=\mathcal{C}_1\odot\mathcal{C}_2$ is a Fra\"{\i}ss\'e class of
finite relational structures. These structures are rigid: for the only finite
group whose order is both odd and a power of $2$ is the trivial group.
Since the structure of any $2$-set is an
arc of the tournament, we see that the automorphism group of the Fra\"{\i}ss\'e
limit is transitive on $2$-element sets. Thus, there cannot be a reduct which
is a total order, since it would agree with the tournament on some but not all
pairs, and the group would not be transitive on $2$-sets. Thus, by the theorem
of Kechris \emph{et al.}, $\mathcal{C}$ is not a Ramsey class.

Here is a direct proof of that fact. First I state the negation of the Ramsey
property: there exist $A$ and $B$ such that, for any $C$, there is a colouring
of the embeddings $A\to C$ red and blue so that there is no monochromatic copy
of $B$. To demonstrate this, we take $A$ to be a $2$-element structure, and
$B$ a $3$-element structure for which the induced tournament is a $3$-cycle.
Given any $C$, we take a total ordering of the points of $C$, and colour the
embeddings $A\to C$ red if the total order agrees with the tournament arc on
the image, blue if it disagrees. Then any $3$-cycle of the tournament on $C$
will have both red and blue embeddings of $A$.

As noted earlier, it is shown in \cite{cl} that there is an explicit failure
of the Ramsey property in any Fra\"{\i}ss\'e class of rigid structures which
does not have a total order as a reduct; we can take $|A|=2$. The following
question is raised:

\begin{qn}
If $\mathcal{C}$ is a Fra\"{\i}ss\'e class of rigid structures which does not
have a total order as a reduct, is there a failure of the Ramsey property with
$|A|=2$ and $|B|$ bounded by a function of the number of isomorphism types of
$2$-element structures?
\end{qn}

For the rest of this section, I will begin a study of Fra\"{\i}ss\'e classes
of rigid structures, which might profitably be taken further. Let $\mathcal{C}$
be such a class and $M$ its Fra\"{\i}ss\'e limit. Then $\Aut(M)$ is a
torsion-free group, since an element of finite order would permute some finite
set non-trivially. From now on I also assume that $M$ has no reduct which is a
total order.

\begin{qn}
Are there examples where $\Aut(M)$ is a simple group?
\end{qn}

Now here are a couple of observations.

\paragraph{Note 1} $M$ has a reduct which is a tournament. For, by homogeneity,
for any $2$-set $\{x,y\}$, there is a relation which orders the pair (else
there would be an automorphism exchanging them); and the collection of all
such relations defines a tournament structure. (Take the relations in turn,
and put $x\to y$ if the first relation to order them satisfies this.)

This tournament is vertex-transitive. We might begin the study by assuming that
the tournament is homogeneous, since all homogeneous tournaments are known
\cite{lachlan}.

\medskip

I will call a ternary relation which distinguishes one of its three arguments
an \emph{AM-relation}, after Samuel Taylor Coleridge's lines~\cite{coleridge}:
\begin{quote}
It is an ancient Mariner,\\
And he stoppeth one of three.
\end{quote}
A C-relation is an example.

\paragraph{Note 2} $M$ has a reduct which is an AM-relation.
For if a triple is totally ordered by the induced tournament on it, we can
take any point to be distinguished. But if the induced tournament is cyclic,
there must be a relation which kills its cyclic automorphism by distinguishing
one point of the three.

\begin{qn}
In our first example, the tournament and AM-relation suffice to make the
structures rigid; presumably this is not true in general. Find examples.
\end{qn}

\begin{qn}
Classify the Fra\"{\i}ss\'e classes with just two relations, a
tournament (which is not a total order) and an AM-relation.
\end{qn}

\end{document}